\title{Elementary abelian regular coverings of Platonic maps\\ Case I: ordinary representations}

\author{Gareth A. Jones\\ School of Mathematics\\ University of Southampton\\ Southampton SO17  1BJ, UK.\\ {\tt G.A.Jones@maths.soton.ac.uk} }

\documentclass[12pt]{article}
\usepackage[latin1]{inputenc}
\usepackage{a4wide}
\usepackage{latexsym}
\usepackage{amsmath}
\usepackage{amssymb}
\usepackage{amsfonts}
\usepackage{color}
\newtheorem{thm}{Theorem}[section]
\newtheorem{lemma}[thm]{Lemma}
\newtheorem{cor}[thm]{Corollary}

\date{}
\begin{document} 

\newcommand{\B}{\mathcal{B}}
\newcommand{\R}{\mathbb{R}}
\newcommand{\Q}{\mathbb{Q}}
\newcommand{\C}{\mathbb{C}}
\newcommand{\Z}{\mathbb{Z}}
\newcommand{\D}{\mathbb{D}}
\newcommand{\E}{\mathfrak{E}}
\newcommand{\F}{\mathbb{F}}
\newcommand{\G}{\mathcal{G}}
\newcommand{\N}{\mathbb{N}}
\newcommand{\K}{\mathcal{K}}
\newcommand{\M}{\mathcal{M}}
\newcommand{\Set}{\mathfrak{Set}}
\newcommand{\Pro}{\mathbb{P}}
\newcommand{\pro}{\mathbb{P}}
\newcommand{\Hyp}{\mathbb{H}}
\newcommand{\hyp}{\mathbb{H}}
\newcommand{\DA}{\mathcal{D}(\mathfrak{A})}
\newcommand{\ug}{\text{\underline{\text{$G$}}}}
\newcommand{\Y}{\mathcal{Y}}
\newcommand{\X}{\mathcal{X}}

\maketitle

\bigskip

\begin{abstract}
We classify the orientably regular maps which are elementary abelian regular branched coverings of Platonic maps $\cal M$, in the case where the covering group and the rotation group $G$ of $\cal M$ have coprime orders. The method involves studying the representations of $G$ on certain homology groups of the sphere, punctured at the branch-points. We give a complete classification for branching over faces (or, dually, vertices) of $\cal M$, and outline how the method extends to other branching patterns.
\end{abstract}




\section{Introduction}

This paper addresses the problem of classifying orientably regular maps, those maps $\cal N$ on surfaces (always assumed to be compact and orientable) for which the orientation-preserving automorphism group ${\rm Aut}^+{\cal N}$ acts transitively on directed edges. Such maps are often formed as regular coverings of a simpler orientably regular map $\cal M$, in which case they correspond to normal subgroups of the fundamental group $\pi_1\,{\cal S}$ of the underlying surface $\cal S$ of $\cal M$, punctured at any branch-points, which are invariant under the induced action of the orientation-preserving automorphism group  $G:={\rm Aut}^+{\cal M}$ on $\pi_1\,{\cal S}$. In the case of abelian coverings, one can abelianise $\pi_1{\cal S}$, and instead look for $G$-invariant submodules of the first homology group $H_1({\cal S};{\mathbb Z})$: for instance, Surowski and the author~\cite{JSu} have used this method to classify the cyclic regular coverings of the Platonic maps, branched over the vertices, edges or faces. When considering coverings by elementary abelian $p$-groups, one can go further and reduce the homology module mod~$(p)$, so that $G$ acts on the vector space $H_1({\cal S};{\mathbb F}_p)$: this idea has been used by Kazaz~\cite{K, K02, K09} to classify the elementary abelian unbranched regular coverings of the orientably regular hypermaps of genus $2$.

Here we combine these two approaches to study regular branched coverings of the Platonic maps $\cal M$, the orientably regular maps of genus $0$, by elementary abelian $p$-groups, as a first step towards a more general classification of the abelian coverings of these maps. The representation theory of the corresponding Platonic groups $G$ plays an important role, and this theory is rather easier to apply in the case where $p$ is coprime to $|G|$, as we shall generally assume here. The remaining cases, where $p$ divides $|G|$ and modular rather than ordinary representation theory is required, will be considered in a separate paper~\cite{J}.

The study of the action of automorphisms on homology is a classical technique in the context of Riemann surfaces~\cite[\S V.3]{FK}, and in recent years it has also been applied to coverings of graphs: see~\cite{KO} and~\cite{MMP}, for instance. In fact, the present work could be restated in terms of coverings of the Platonic graphs, restricted to those coverings which respect the surface-embeddings.

The main result is the following theorem, where we use the notation $\{n,m\}$ of Coxeter and Moser~\cite{CM} for the $m$-valent Platonic map with $n$-gonal faces:

\begin{thm} The orientably regular coverings of the Platonic maps $\cal M$, branched over the faces, with elementary abelian $p$-groups of rank $c\ge 1$ as covering groups, where the prime $p$ does not divide the order of the group $G={\rm Aut}^+{\cal M}$, are as follows:
\begin{itemize}
\item the tetrahedron $\{3,3\}$ has one regular covering for each $p$, with $c=3$;
\item the cube $\{4,3\}$ has three regular coverings for each $p$, with $c=2,3,5$;
\item the octahedron $\{3,4\}$ has seven regular coverings for each $p$, with $c=1,3,3,4,4,6,7$;
\item the dodecahedron $\{5,3\}$ has seven regular coverings for each $p\equiv\pm 1$ {\rm mod}~$(5)$, with $c=3,3,5,6,8,8,11$, and has three regular coverings for each $p\equiv\pm 2$ {\rm mod}~$(5)$, with $c=5,6,11$;
\item the icosahedron $\{3,5\}$ has $8p+23$ coverings, $31$ regular and $8(p-1)$ chiral, for each $p\equiv\pm 1$ {\rm mod}~$(5)$, with $c\in\{3,4, \ldots, 16, 19\}$, and has $4p+11$ coverings, $15$ regular and $4(p-1)$ chiral, with $c\in\{4,5,6,8,9,10,11,13,14,15,19\}$, for each $p\equiv\pm 2$ {\rm mod}~$(5)$;
\item the dihedron $\{n,2\}$ has one regular covering for each $p$, with $c=1$;
\item the hosohedron $\{2,n\}$ has $2^{\nu}-1$ regular coverings for each $p$, where $\nu$ is the number of orbits $\Delta\ne\{1\}$ of the group generated by the Frobenius automorphism and inversion on the set of $n$th roots of $1$ in the algebraic closure $\overline{\mathbb F}_p$ of the field ${\mathbb F}_p$.
\end{itemize}
For each ${\cal M}=\{n,m\}$ and prime $p$ these covering maps have type $\{np,m\}$ and genus
\[\left(\frac{f}{2}-1\right)p^c-\frac{f}{2}p^{c-1}+1,\]
where $f$ is the number of faces of $\cal M$; they are all quotients of a single covering with $c=f-1$. In all cases except when $\cal M$ is the icosahedron, the $2^{\nu}-1$ coverings are the joins formed from a set of $\nu$ irreducible coverings of $\cal M$.
\end{thm}

\noindent{\bf Comments. 1.} The Schur-Zassenhaus Theorem implies that the orientation-preserving automorphism group of each of these coverings is a semidirect product of the covering group, an elementary abelian normal subgroup of order $p^c$, by $G$. 

\smallskip

\noindent{\bf 2.} In the case of the hosohedron, the method for finding the ranks $c$ of the covering groups, too complicated to state here, is described in \S 9.2.

\smallskip

\noindent{\bf 3.} A map is called regular or chiral as it is or is not isomorphic to its mirror image. When $\cal M$ is the icosahedron, the chiral coverings arise from the fact that its full automorphism group ${\rm Aut}\,{\cal M}$ has fewer orbits than $G$ on pairs of faces (see \S 2.6 and \S 7). The fact that the number of coverings is unbounded as a function of $p$ is related to the fact that the permutation character of $G$ on faces is not multiplicity-free (see \S 2.5 and \S 7).

\smallskip

\noindent{\bf 4.} The conclusions of Theorem~1.1 also apply in a few cases where $p$ divides $|G|$, such as the tetrahedron and octahedron with $p=3$, and the hosohedron with $n$ odd and $p=2$.

\medskip

The paper is organised as follows. In \S 2 some useful techniques are outlined. These are applied in \S\S3--9 to enumerate and describe the coverings in Theorem~1.1. The individual cases are fairly straightforward, except for that involving the hosohedron ${\cal M}=\{n,2\}$, which depends on the reduction mod~$(p)$ of cyclotomic polynomials. Branching over vertices is easily dealt with by duality. In \S 10 we consider simultaneous branching over vertices and faces, concentrating on the tetrahedron as a typical example. In \S\S11--12 we briefly consider the hypermaps arising from branching over edges (and also vertices and faces).

The author is grateful to Young Soo Kwon for raising this problem, and to Rub\'en Hidalgo and Roman Nedela for helpful comments on an earlier draft of this work.

\section{Techniques}

\subsection{Maps}

First we briefly sketch the connections between maps (always assumed to be orientable) and triangle groups; for background, see~\cite{JSi}. We define a triangle group to be
\[\Delta(p, q, r)=\langle x, y, z \mid x^p=y^q=z^r=xyz=1 \rangle\]
where $p,q,r\in{\N}\cup\{\infty\}$ and we ignore any relation $g^{\infty}=1$. Any $m$-valent map $\cal N$ corresponds to a subgroup $N$ of the triangle group
\[\Delta:=\Delta(m, 2, \infty)=\langle x, y, z \mid x^m=y^2=xyz=1\rangle,\]
with vertices, edges and faces corresponding to the cycles of $x, y$ and $z$ on the cosets of $N$. The map $\cal N$ is orientably regular if and only if $N$ is normal in $\Delta$, in which case the orientation-preserving automorphism group ${\rm Aut}^+{\cal N}$ is isomorphic to $\Delta/N$; we will assume this throughout. In particular, the Platonic map ${\cal M}=\{n, m\}$ corresponds to the normal closure $M$ of $z^n$ in $\Delta$, and ${\rm Aut}^+{\cal M}\cong\Delta/M\cong\Delta(m,2,n)$. A map $\cal N$ is a $d$-sheeted covering of $\cal M$ if and only if $N$ is a subgroup of index $d$ in $M$, in which case ${\cal N}\to{\cal M}$ is a regular covering and the group of covering transformations is $M/N$. The orientably regular $m$-valent maps which are $d$-sheeted abelian coverings of $\cal M$ are therefore in bijective correspondence with the subgroups $N$ of $M$ which are normal in $\Delta$, with $M/N$ abelian of order $d$. If $M/N$ has exponent $l$, such subgroups $N$ contain the commutator subgroup $M'$ and the subgroup $M^l$ generated by the $l$-th powers in $M$, so they correspond to subgroups $\overline N=N/M'M^l$ of $\overline M=M/M'M^l$, with $M/N\cong \overline M/\overline N$. The action of $\Delta$ by conjugation on the normal subgroup $M$ preserves its characteristic subgroups $M'$ and $M^l$, so there is an induced action of $\Delta$ on $\overline M$; since $M$ is in the kernel of this action, we therefore have an action of the group $G:={\rm Aut}^+{\cal M}\cong\Delta/M$ on $\overline M$, which is a module for $G$ over the ring ${\Z}_l$; it follows that a subgroup $N$ of $M$, containing $M'M^l$, is normal in $\Delta$ if and only if $\overline N$ is a $G$-invariant submodule of $\overline M$. The orientation-preserving automorphism group $\tilde G:={\rm Aut}^+{\cal N}\cong\Delta/N$ of $\cal N$ has an abelian normal subgroup $K$, corresponding to $M/N$, with $\tilde G/K\cong G$; under the induced action of $G$ on $K$ by conjugation, $K$ is isomorphic as a $G$-module to $\overline M/\overline N$. If $l$ is coprime to $|G|$ then the Schur-Zassenhaus Theorem implies that $\tilde G$ is a semidirect product of $K$ by $G$.

\subsection{Homology modules}

Let $S=S^2\setminus\{c_1,\ldots, c_f\}$, where $c_1,\ldots, c_f$ are the centres of the $f$ faces of $\cal M$. Then $M$ can be identified with the fundamental group $\pi_1(S)$ of $S$, a free group of rank $f-1$ generated by the homotopy classes $g_i$ of loops around the punctures $c_i$, with a single defining relation $g_1\ldots g_f=1$. It follows that the group $M^{\rm ab}=M/M'$ can be identified with the first integer homology group $H_1(S;{\Z})=\pi_1(S)^{\rm ab}$ of $S$, a free abelian group of rank $f-1$ generated by the homology classes $[g_i]$, with defining relation $[g_1]+\cdots+[g_f]=0$. By the Universal Coefficient Theorem, $\overline M$ is identified with the mod~$(l)$ homology group $H_1(S;{\Z}_l)=H_1(S;{\Z})\otimes_{\mathbb Z}{\Z}_l\cong{\Z}_l^{f-1}$. Under these identifications, the actions of $G$ induced by conjugation in $\Delta$ and by homeomorphisms of $S$ are the same, so our problem is to find the $G$-submodules of $H_1(S;{\Z}_l)$.

The prime power factorisation $\prod_pp^{e_p}$ of $l$ induces a $G$-invariant decomposition of $H_1(S;{\Z}_l)$ as a direct sum of its Sylow $p$-subgroups $H_1(S;{\Z}_{p^{e_p}})$, so it is sufficient to restrict attention to the case where $l$ is a power of a prime $p$. For each prime $p$ we have an infinite descending series
\[M^{\rm ab} = M/M' > M'M^p/M' > M'M^{p^2}/M' > \ldots > M'M^{p^i}/M' > \ldots,\]
of characteristic subgroups of $M^{\rm ab}$, corresponding to a descending series
\[H_1(S;{\Z}) = H^0 > H^1>\ldots > H^i > \ldots,\]
of $G$-submodules of $H_1(S;{\Z})$, where each $H^i$ is the kernel of the reduction mod~$(p^i): H_1(S;{\Z})\to H_1(S;{\Z}_{p^i})$. This induces a finite descending series
\[H_1(S;{\Z}_{p^i}) = H^0/H^i > H^1/H^i>\ldots > H^i/H^i = 0\]
of $G$-submodules of $H_1(S;{\Z}_{p^i})$. Successive quotients
\[(H^j/H^i)/(H^{j+1}/H^i) \cong H^j/H^{j+1} \cong M'M^{p^j}/M'M^{p^{j+1}}\]
in this series are $G$-isomorphic to the module
\[H_1(S;{\F}_p) = H_1(S;{\Z}_p) = H^0/H^1 \cong M_p:= M/M'M^p,\]
where $\F_p$ is the field of $p$ elements.
This $G$-module $H_1(S;{\F}_p)$, along with all the other quotients in the series, is a vector space over $\F_p$ of dimension $f-1$ affording a representation $\rho_p:G\to {\rm Aut}(H_1(S;{\F}_p)) \cong GL_{f-1}(\F_p)$ of $G$.

Our main problem therefore is to determine the structure of this $G$-module $H_1(S;{\F}_p)$ for each prime $p$. In particular, this will immediately determine the elementary abelian regular coverings of $\M$. A $G$-submodule of $L$ of $H_1(S;{\F}_p)$ of codimension $c$ corresponds to a normal subgroup $N$ of $\Delta$ contained in $M$, with $M/N$ elementary abelian of order $p^c$, and hence to a regular $p^c$-sheeted covering $\cal N$ of $\cal M$, branched over its face-centres. This is an orientably regular map of type $\{np, m\}$. There are $p^{c-1}$ points above each of the $f$ face-centres of $\cal M$, so the total order of branching of the covering ${\cal N}\to{\cal M}$ is $(p^c-p^{c-1})f$. By the Riemann-Hurwitz formula, $\cal N$ therefore has genus
\[g=1-p^c+\frac{1}{2}p^{c-1}(p-1)f= \left(\frac{f}{2}-1\right)p^c-\frac{f}{2}p^{c-1}+1.\]
The orientation-preserving automorphism group $\tilde G={\rm Aut}^+{\cal N}$ of $\cal N$ has an elementary abelian normal subgroup $K$ of order $p^c$, isomorphic to $H_1(S;{\F}_p)/L$ as a $G$-module over $\F_p$.

For notational convenience, we let $E_p(\M)$ denote the set of such regular coverings of $\M$ with  elementary abelian $p$-groups as covering groups, and for each ${\cal N}\in E_p(\M)$ we let $\dim({\cal N})=c$ where the covering ${\cal N}\to\M$ has $p^c$ sheets. Thus $\dim({\cal N})$ is simply the dimension of the corresponding $G$-module $H_1(S;\F_p)/L$, or equivalently the codimension of $L$. We will denote the representations of $G$ on $L$ ($\,\cong_GN/M'M^p$) and on $H_1(S;{\F}_p)/L$ ($\,\cong_G M/N \cong_G K$) respectively by $\rho_L$ and $\rho^L$, and the corresponding characters by $\chi_L$ and $\chi^L$. In particular, we will say that the covering $\cal N$ {\sl affords\/} the representation $\rho^L$ and the character $\chi^L$, since these correspond to the action of $G$ by conjugation on the covering group $K$. We will call a covering $\cal N$ {\sl irreducible\/} or {\sl indecomposable\/} if $\rho^L$ has this property.

 If $p$ is coprime to $|G|$ then the representation theory of $G$ over fields of characteristic $p$ is essentially the same as that over fields of characteristic $0$ (ordinary representation theory): for instance Maschke's Theorem applies in both cases. In this situation, the decomposition of $H_1(S;\F_p)$ as a $G$-module can be obtained from that of the corresponding homomology module over $\mathbb C$ (or indeed any algebraically closed field of characteristic $0$, such as the field $\overline{\mathbb Q}$ of algebraic numbers). Moreover, if $p$ is coprime to $|G|$ then $\tilde G$ is a semidirect product of $K$ by $G$, by the Schur-Zassenhaus Theorem.
 
 If, on the other hand, $p$ divides $|G|$ then this extension splits if and only if its Sylow $p$-subgroups split over $K$. Moreover, in this situation, Maschke's Theorem does not apply, and the submodule structure of $H_1(S;\F_p)$ is less transparent, requiring modular rather than ordinary representation theory. In this paper we therefore concentrate on the ordinary case, avoiding the finitely many primes $p$ dividing $|G|$; these will be dealt with in a later paper~\cite{J}.
 
 \subsection{Permutation modules}

Let $\Phi$ be the set of faces of $\cal M$, and let $\C\Phi$ be the corresponding permutation module for $G$; this is an $f$-dimensional complex vector space with basis $\Phi$ permuted naturally by $G$. The homology module $H_1(S;{\mathbb C})$ is isomorphic to the quotient of $\C\Phi$ by the $1$-dimensional $G$-submodule spanned by the element $\sum_{\phi\in\Phi}\phi$. The $G$-module decomposition of $\C\Phi$ corresponds to that of the corresponding permutation character $\pi$, where $\pi(g)$ is the number of faces $\phi\in\Phi$ invariant under each $g\in G$. If $H$ denotes the subgroup $\langle z\rangle\cong C_n$ of $G$ leaving invariant a face, then since $G$ acts transitively on $\Phi$ we have $\pi=1_H^G$, the character obtained by inducing the principal character $1_H$ of $H$ up to $G$. By Frobenius reciprocity, the multiplicity $(\chi_i, \pi)$ of any irreducible complex character $\chi_i$ of $G$ in $\pi$ is equal to $(\chi_i\mid_H, 1_H)$, the multiplicity of $1_H$ in the restriction of $\chi_i$ to $H$. This is equal to the average value $|H|^{-1}\sum_{g\in H}\chi_i(g)$ of $\chi_i$ on $H$, or equivalently the multiplicity of $1$ as an eigenvalue for $z$ in the irreducible representation $\rho_i$ of $G$ corresponding to $\chi_i$. Once these multiplicities are calculated, the direct sum decomposition of  $\C\Phi$ is known, and hence so is that of $H_1(S;{\mathbb C})$. In particular, this module affords the character
\[\chi=\pi-\chi_1\]
of $G$, where $\chi_1$ is the principal character $1_G$ of $G$, given by $\chi_1(g)=1$ for all $g\in G$.

The same decompositions apply if $\C$ is replaced with the algebraic closure $\overline{\F}_p$ of $\F_p$, where $p$ is any prime not dividing $|G|$. In order to find the corresponding decompositions over $\F_p$, various algebraically conjugate summands must be merged to give summands which are realised over $\F_p$. For any given map Platonic $\M$, how this happens depends on certain congruences satisfied by $p$.


\subsection{Finding submodules}

Some submodules $L$ of $H_1(S;\F_p)$ arise naturally from the action of $G$ on $\Phi$, whether or not $p$ divides $|G|$. Let $P$ denote the permutation module $\F_p\Phi$, a $G$-module over $\F_p$. Given any subset $\Psi\subseteq\Phi$, let $\underline\Psi$ denote the element $\sum_{\phi\in\Psi}\phi$ of $P$. The $G$-module $H_1(S;\F_p)$ can then be identified with the quotient $Q:=P/P_1$ of $P$ by the $1$-dimensional $G$-submodule $P_1$ spanned by $\underline\Phi$.

There is a regular covering $\M_0\in E_p(\M)$, with $\dim({\M_0})=f-1$, corresponding to the submodule $L=0$ of $Q=H_1(S;\F_p)$, i.e.~to the submodule $P_1$ of $P$; all other regular coverings $\cal N$ of $\cal M$ in $E_p(\M)$ are proper quotients of $\M_0$, so they satisfy $\dim({\cal N})\le f-2$.

There is a $G$-submodule $P^1$ of codimension $1$ in $P$, consisting of the elements $\sum_{\phi\in\Phi}a_{\phi}\phi$ with coordinate-sum $\sum_{\phi\in\Phi}a_{\phi} =0$. If $p$ divides $f$ then $P_1\le P^1$, giving a $G$-submodule $L=Q^1=P^1/P_1$ of codimension $1$ in $Q$, and hence a regular covering $\M^1\in E_p(\M)$ with $\dim(\M^1)=1$. This is therefore a cyclic covering of $\M$, and since $L/Q^1\cong P/P^1$ affords the principal representation of $G$ it is also a central covering, in the sense that $K$ is in the centre of $\tilde G$. The cyclic coverings of the Platonic hypermaps were classified by the author and Surowski in~\cite{JSu} (see also~\cite{SJ}), where Proposition~3 implies that $\tilde G$ is a split extension of $G$ by $K$, giving $\tilde G=G\times K\cong G\times C_p$, if and only if $p$ divides $2m$ but not $n$. On the other hand, if $p$ does not divide $f$ then $P=P_1\oplus P^1$ so that $Q\cong P^1$ and we obtain the identity covering $\M\to\M$.

If $G$ acts imprimitively on $\Phi$, preserving a non-trivial equivalence relation $\sim$ (equivalently, if there is a subgroup $H_{\sim}$ of $G$ such that $H<H_{\sim}<G$), then as $\Psi$ ranges over the $k=|G:H_{\sim}|$ equivalence classes, the elements $\underline{\Psi}$ form a basis for a $k$-dimensional $G$-submodule $P_{\sim}$ of $P$, containing $P_1$ (and contained in $P^1$ if and only if $p$ divides the size $|H_{\sim}:H|=f/k$ of each class). The $G$-submodule $L=Q_{\sim}=P_{\sim}/P_1$ of $Q$ therefore corresponds to a regular covering $\M_{\sim}\in E_p(\M)$ with $\dim(\M_{\sim})=f-k$. The normal subgroup $K$ of $\tilde G$ affords a representation of $G$ with character $\pi-\pi_{\sim}$, where $\pi_{\sim}$ is the permutation character of $G$ on the equivalence classes of $\sim$, i.e.~on the cosets of $H_{\sim}$ in $G$. If $\approx$ is a $G$-invariant equivalence relation which refines $\sim$ (i.e.~$\phi_1\approx\phi_2$ implies $\phi_1\sim\phi_2$), then the inclusion $Q_{\sim}\le Q_{\approx}$ induces a covering $\M_{\sim}\to\M_{\approx}$.

For instance, suppose that $\cal M$ has an antipodal symmetry, so that $f$ is even. Then $P$ has a $G$-submodule $P_a\ge P_1$ of dimension $k=f/2$ with basis elements $\underline\Psi$ corresponding to the antipodal pairs $\Psi\subseteq\Phi$. This gives a $G$-submodule $L=Q_a=P_a/P_1$ of codimension $f-k=f/2$ in $Q$, so we obtain a regular covering $\M_a\in E_p(\M)$ with $\dim(\M_a)=f/2$. If $p>2$ there is a second $G$-submodule $P_{a'}$  of dimension $f/2$ in $P$, with a basis element $\phi-\phi'$ for each antipodal pair $\{\phi, \phi'\}\subseteq\Phi$ (so $P_{a'}\le P^1$). The $G$-submodule $L=Q_{a'}=(P_{a'}\oplus P_1)/P_1$ of codimension $(f-2)/2$ in $Q$ corresponds to a regular covering $\M_{a'}\in E_p(\M)$ with $\dim(\M_{a'})=(f-2)/2$. A simple calculation shows that $P=P_a\oplus P_{a'}$, so $Q=Q_a\oplus Q_{a'}$. Since $P_a$ affords the permutation character $\pi_a$ of $G$ on the set of antipodal pairs of faces, it follows that $Q_a$ and hence $Q/Q_{a'}$ afford $\pi_a-\chi_1$, while $Q_{a'}$ and $Q/Q_a$ afford $\chi-\pi_a$.



If $p=2$ one can identify $P$ with the power set ${\cal P}(\Phi)$ of $\Phi$, a group under symmetric difference, with the natural induced action of $G$, so that $P^1=\{\Psi\subseteq\Phi\mid|\Psi|\;\hbox{is even}\}$ and $P_1=\{\emptyset, \Phi\}$. In this case, $Q$ can be identified with the set of complementary pairs $\{\Psi, \Phi\setminus\Psi\}$ of subsets of $\Phi$.

For each $\M$, the group $G$ has a natural $3$-dimensional complex representation $\rho_n$, obtained by extending its natural real representation, as the rotation group of $\M$, by linearity to $\C$. This representation is irreducible, except when $\M$ is the dihedron or hosohedron, in which case it splits as a direct sum of $1$- and $2$-dimensional representations. When $\M$ is one of the five Platonic solids, by sending each basis element $\underline\phi\;(\phi\in\Phi)$ of $\C\Phi$ to the position-vector in $\R^3\subset\C^3$ of the face-centre of $\phi$, we see that $\rho_n$ is a quotient, and hence a direct summand, of the representation of $G$ on $\C\Phi$. Frobenius reciprocity implies that $\rho_n$ has multiplicity $1$ in this representation, and in the cases of the dodecahedron and icosahedron, the same applies to the Galois conjugate representation $\rho_{n*}$. By realising $\rho_n$ over the ring of integers of an appropriate algebraic number field (namely $\Q$ for the tetrahedron, cube or octahedron, $\Q(\sqrt 5)$ for the dodecahedron or icosahedron, and a cyclotomic field $\Q(\zeta_n)$ for the dihedron or hosohedron, where $\zeta_n=e^{2\pi i/n}$), and then taking quotients modulo some prime ideal containing $p$, we find that the mod~$(p)$ reduction of $\rho_n$ or of $\rho_n\oplus\rho_{n*}$ is a summand of the representation of $G$ on $P$, and hence on $Q$.

\subsection{Diagonal submodules}

We need to deal with those cases when $Q$ has summands of multiplicity greater than one. If $A_1$ and $A_2$ are isomorphic $FG$-modules for some field $F$ and group $G$, then a {\sl diagonal submodule\/} of $A_1\oplus A_2$ is a submodule of the form $\{(a,a')\mid a\in A_1\}$ for some isomorphism $A_1\to A_2, a\mapsto a'$.

\begin{lemma}
Let $U$ be a submodule of a module $V=V_1\oplus V_2$, let $U_i=U\cap V_i$, and let $\overline V_i=V_i/U_i$ for $i=1, 2$. Then $U$ is the inverse image in $V$ of a diagonal submodule of $\overline W_1\oplus\overline W_2$ for some submodules $W_i$ such that $U_i\le W_i\le V_i$ for $i=1, 2$.
\end{lemma}

\noindent{\sl Proof.} The image of $U$ in $\overline V_1\oplus\overline V_2$ intersects each summand $V_i$ trivially, so it is isomorphic to its projection $\overline W_i$ in each $V_i$, and is therefore a diagonal submodule of $\overline W_1\oplus\overline W_2$. \hfill$\square$

\medskip

Thus, provided the submodule structures of each $V_i$ are known, so is that of $V$.

\begin{cor}
If $V_1$ and $V_2$ are isomorphic irreducible $FG$-modules, then the only proper submodules of $V_1\oplus V_2$ are $0$, $V_1$, $V_2$ and the diagonal submodules; given any isomorphism $V_1\to V_2, v\mapsto v'$, each of the latter has the form $\{(v,\lambda v')\mid v\in V_1\}$ for some $\lambda\in F^*$.
\end{cor}

\noindent{\sl Proof.} Let $U$ be a submodule of $V=V_1\oplus V_2$. In the notation of Lemma~1.1, for each $i=1, 2$ the irreducibility of $V_i$ implies that $U_i=0$ or $V_i$. If $U_i=V_i$ for some $i$ then $U\ge V_i$, so the irreducibility of $V/V_i$ implies that $U=V_i$ or $V$. Otherwise, $U_1=U_2=0$, so in Lemma~1.1 we have either $W_i=0$ for some $i$, giving $U=0$, or $W_i=V_i$ for each $i$, in which case $U$ is a diagonal submodule of $V$. In the latter case, Schur's Lemma  implies that the isomorphisms $V_1\to V_2$ are the functions $v\mapsto\lambda v'$ for $\lambda\in F^*$, so $U$ is as claimed.\hfill$\square$

\medskip

This shows that if $V_1=V_2$ is irreducible and $F=\F_q$, then $V=V_1\oplus V_2=V_1\oplus V_1$ has $q+1$ non-trivial proper submodules, all isomorphic to $V_1$. We can write them as
\[V(\lambda)=\{(\alpha v, \beta v)\mid v\in V_1,\, \alpha/\beta=\lambda\}\]
for $\lambda\in{\mathbb P}^1(\F_q)=\F_q\cup\{\infty\}$; here $V(\infty)$ and $V(0)$ are the first and second direct summands, with $\beta=0$ and $\alpha=0$ respectively.

\subsection{Regularity and chirality}

Each of the coverings $\cal N$ constructed here is orientably regular, since it corresponds to a normal subgroup $N$ of $\Delta$. It is regular (i.e.~it also admits an orientation-reversing automorphism) if and only if $N$ is normal in the extended triangle group $\Gamma=\Delta[m, 2, \infty]$, which contains $\Delta$ as a subgroup of index $2$. This is equivalent to the corresponding submodule $L$ of $H_1(S;\Z_l)$ being invariant, not just under the orientation-preserving automorphism group $G={\rm Aut}^+\M\cong \Delta/N$ of $\M$, but under its full automorphism group $A={\rm Aut}\,\M\cong\Gamma/N$. In the case $l=p$, it is clear from their construction that the coverings $\M_0$, $\M^1$, $\M_a$ and $\M_{a'}$ (when they exist) are all regular.

One can find the structure of $H_1(S;\F_p)$ as an $A$-module by treating it as the quotient $Q=P/P_1$ of $P$, where $P$ is now regarded as the permutation module for $A$ on $\Phi$. In particular, the permutation character $\pi|_A$ of $A$ on $P$ can be found as for $\pi|_G$, taking $H$ to be the stabiliser $A_{\phi}\cong D_n$ of a face $\phi\in\Phi$ in $A$ instead of its stabiliser $G_{\phi}\cong C_n$ in $G$.

If $\M$ is not the icosahedron, and if $p$ is a prime not dividing $|G|$, then all the coverings ${\cal N}\in E_p(\M)$ are regular. To prove this, it is sufficient to show that the $G$-submodules of $H_1(S;\C)$ are all $A$-invariant, i.e.~that the $A$- and $G$-module direct sum decompositions of $H_1(S;\C)$ are identical. If this is not the case, then at least one irreducible summand of this $A$-module must split into a direct sum of two irreducible summands of the $G$-module, so that the permutation character $\pi$ on $\Phi$ satisfies $(\pi|_A,\pi|_A)_A<(\pi|_G,\pi|_G)_G$. But the two sides of this inequality are equal to the ranks of $A$ and $G$ as permutation groups on $\Phi$, i.e.~the number of orbits of $A$ and $G$ on $\Phi^2$, or equivalently the number of orbits of their face-stabilisers $A_{\phi}$ and $G_{\phi}$ on $\Phi$. Now $A_{\phi}\cong D_n$ and $G_{\phi}\cong C_n$, and for all $\M$ except the icosahedron their orbits on $\Phi$ are either trivial (fixing $\phi$ and its antipodal face, if it exists), or natural orbits of length $n$. Thus $A_{\phi}$ and $G_{\phi}$ have the same number of orbits on $\Phi$, so $A$ and $G$ have the same rank, and the assertion is proved.
This argument fails for the icosahedron, with $A$ and $G$ having ranks $6$ and $8$ on $\Phi$, and $A_{\phi}$ having two regular orbits, of length $2n=6$; later we will give explicit examples of chiral maps in $E_p(\M)$ in this case.

\section{The tetrahedron}\label{Tet}

Let $\M$ be the tetrahedral map $\{3, 3\}$, with $f=4$ faces and $G={\rm Aut}^+\M\cong A_4$. The group $A_4$ has four conjugacy classes, consisting of the identity, the three double transpositions, and two mutually inverse classes of four $3$-cycles. Its character table is as follows; here $\chi_1$ is the principle character, $\chi_2$ and $\chi_3$ are the faithful characters of $A_4/V_4\cong C_3$ where $V_4$ is the normal Klein four-group, and $\chi_4$ is the character of the natural representation $\rho_n$ as the rotation group of the tetrahedron, extended from $\R^3$ to $\C^3$.

\begin{table}[htb]
\centering
\begin{tabular}{|c|c|c|c|c|}
\hline
&$1$&$(..)(..)$&$(...)^+$&$(...)^-$\\
\hline
$\chi_1$&1&1&1&1\\
$\chi_2$&1&1&$\omega$&$\overline\omega$\\
$\chi_3$&1&1&$\overline\omega$&$\omega$\\
$\chi_4$&3&-1&0&0\\ \hline
\end{tabular}
\caption{The character table of $A_4$, with $\omega=\exp(2\pi i/3)$.}\label{chartA4}
\end{table}

Since $G$ acts doubly transitively on $\Phi$, as the natural representation of $A_4$, it follows that $\pi=\chi_1+\chi_4$, so $H_1(S;\C)$ is an irreducible $G$-module affording the character $\chi_4$.

The primes dividing $|G|$ are $2$ and $3$, so for primes $p>3$ we can regard Table~\ref{chartA4} as the character table of $A_4$ over $\overline\F_p$ by suitable reduction mod~$(p)$, with $\omega^2+\omega+1=0$. In particular, the $G$-module $Q=H_1(S;\F_p)$ is irreducible, affording the character $\chi_4$ of $G$. It follows that in this case $E_p(\M)$ consists of a single regular map $\M_0$, with $\dim(\M_0)=3$. (If $p=3$ then $Q$ is again irreducible, and the result is the same.) The map $\M_0$ has type $\{3p, 3\}$ and genus $p^3-2p^2+1$.
Its orientation-preserving automorphism group ${\rm Aut}^+\M_0$ has an elementary abelian normal subgroup $K\cong (C_p)^3$; the quotient group, isomorphic to $G$, has an induced action by conjugation on $K$, so that $K$ is a $G$-module over $\F_p$, isomorphic to $Q$. The full automorphism group ${\rm Aut}\,\M_0$ is an extension of $K$ by $A={\rm Aut}\,\M\cong S_4$, with $K\cong Q=P/P_1$ where we now regard $P$ as the natural permutation module for $S_4$; the character of $S_4$ on $K$ is that denoted by $\chi_4$ in Table~2 (see \S 4). When $p=3$ or $5$ the map $\M_0$ has genus $10$ or $76$, and appears as the dual of the map R10.1 or R76.2 in Conder's computer-generated list of regular orientable maps~\cite{C}.

{\color{red}
}


\section{The cube}\label{Cub}

Next we consider coverings of the cube ${\cal M}=\{4, 3\}$ branched over its face-centres (or, dually, of the octahedron $\{3, 4\}$, branched over its vertices). In this case $G=\Delta(3, 2, 4)\cong S_4$, acting on the set $\Phi$ of six faces of $\cal M$, so the stabilisers of faces are the three subgroups $H\cong C_4$ of $G$.

In $S_4$ there are five conjugacy classes: the identity, six transpositions, three double transpositions, eight $3$-cycles, and six $4$-cycles. Hence there are five irreducible characters. In addition to the principal character $\chi_1$ and the alternating character $\chi_2(g)={\rm sgn}(g)$,  $\chi_3$ and $\chi_4$ are the non-principal irreducible characters obtained from the doubly transitive representations of $S_4$ of degrees $3$ (via the epimorphism $S_4\to S_3$) and $4$ (the natural action), while $\chi_5$ is the character of the natural representation $\rho_n$ of $G$ as the rotation group of the cube or octahedron.

\begin{table}[htb]
\centering
\begin{tabular}{|c|c|c|c|c|c|}
\hline
&$1$&$(..)$&$(..)(..)$&$(...)$&$(....)$\\
\hline
$\chi_1$&1&1&1&1&1\\
$\chi_2$&1&-1&1&1&-1\\
$\chi_3$&2&0&2&-1&0\\
$\chi_4$&3&1&-1&0&-1\\
$\chi_5$&3&-1&-1&0&1\\ \hline
\end{tabular}
\caption{The character table of
$S_4$.}\label{chartS4}
\end{table}

The subgroup $H$ consists of the identity element, one double transposition and two $4$-cycles. By averaging their values over $H$, one therefore sees that the characters $\chi_i$ have multiplicities $1, 0, 1, 0$ and $1$ in the permutation character $\pi$ on $P$, that is,
\[\pi=\chi_1+\chi_3+\chi_5.\]
It follows that the representation of $G$ on $H_1(S;{\mathbb C})$ has character
\[\chi=\chi_3+\chi_5,\]
so $H_1(S;{\mathbb C})$ is a direct sum of irreducible submodules of dimensions $2$ and $3$.

The primes dividing $|G|$ are $2$ and $3$, so we can use this decomposition to find the decomposition of $H_1(S;{\mathbb F}_p)$ for any prime $p>3$. In this case $H_1(S;{\mathbb F}_p)$ is also a direct sum of irreducible submodules of dimensions $2$ and $3$, so $E_p(\M)$ consists of three covering maps $\cal N$, with $\dim({\cal N})=2$, $3$ and $5$. Since the cube has an antipodal symmetry, with the $f=6$ faces permuted in three antipodal pairs, the first two coverings $\cal N$ are the maps $\M_{a'}$ and $\M_a$ corresponding to the summands $L=Q_{a'}$ and $Q_a$ in the $G$-module decomposition $Q=Q_a\oplus Q_{a'}$, while the third is $\M_0$, corresponding to $L=0$. These are regular orientable maps of type $\{4p, 3\}$ and genus $2p^c-3p^{c-1}+1$ where $c=\dim({\cal N})=2, 3$ or $5$, affording the character $\chi_3$, $\chi_5$ or $\chi_3+\chi_5$ respectively; for instance, if $p=5$ or $7$ then the maps $\M_{a'}$ are the duals of the regular maps R36.3 and R78.1 in~\cite{C}. The inclusions $0<Q_{a'}$ and $0<Q_a$ of $G$-submodules induce regular coverings $\M_0\to\M_{a'}$ and $\M_0\to M_a$, with covering groups $Q_{a'}$ and $Q_a$.


{\color{red}
}

{\color{red}
}

\section{The octahedron}\label{Oct}

We now consider coverings of the octahedron ${\cal M}=\{3, 4\}$ branched over its face-centres (or, dually, of the cube $\{4, 3\}$, branched over its vertices). As in the case of the cube we have $G\cong S_4$, but now acting on the set $\Phi$ of eight faces of $\cal M$, so the stabilisers of faces are the four Sylow $3$-subgroups $H\cong C_3$ of $G$. The conjugacy classes and character table of $G$ are as before, but now $H$ consists of the identity and two $3$-cycles. We find that
\[\pi=\chi_1+\chi_2+\chi_4+\chi_5,\]
so
\[\chi=\chi_2+\chi_4+\chi_5.\]
Thus $H_1(S;{\mathbb C})$ is a direct sum of irreducible $G$-submodules of dimensions $1, 3$ and $3$. For primes $p>3$ we therefore obtain seven coverings ${\cal N}\in E_p(\M)$ with $c=\dim({\cal N})=1, 3, 3, 4, 4, 6$ and $7$, corresponding to the seven proper $G$-submodules $L$ of $Q$. These are regular maps of type $\{3p, 4\}$ and genus $g=3p^c-4p^{c-1}+1$. When $c=1$ we have a regular map $\cal N$ of genus $3(p-1)$, affording the character $\chi_2$ of $G$. The underlying Riemann surfaces of these maps are examples of the well-known Accola-Maclachlan surfaces~\cite{Acc, Macl} with $8(g+3)$ automorphisms; for instance, taking $p=5, 7, 11, 13, 17, 19, 23, 29, 31$ we obtain the duals of the regular maps R12.1, R18.1, R30.1, R36.5, R48.1, R54.1, R66.3, R84.1, R90.1 in~\cite{C}. In each of the cases $c=3$ or $c=4$ we obtain a pair of $p^c$-sheeted coverings $\cal N$ of $\cal M$; not only are they non-isomorphic, but their automorphism groups are non-isomorphic since the corresponding pairs of normal subgroups $K$ of $\tilde G$ are non-isomorphic as $G$-modules, affording characters $\chi_4$ and $\chi_5$ when $c=3$, or $\chi_2+\chi_4$ and $\chi_2+\chi_5$ when $c=4$. The remaining coverings, with $c=6$ and $7$, afford the characters $\chi_4+\chi_5$ and $\chi_2+\chi_4+\chi_5$.

These $G$-submodules $L\le Q$, and the corresponding maps ${\cal N}\in E_p({\cal M})$, all have combinatorial interpretations. The octahedron has an antipodal symmetry, so for each prime $p$ there is a $G$-submodule $Q_a$ affording the character $\pi_a-\chi_1$ ($=\chi_4$ if $p>3$), giving a regular map ${\cal N}=M_a\in E_p(\M)$ with $\dim({\cal N})=4$ and $K\cong Q/Q_a$ affording $\chi-\pi_a+\chi_1$ ($=\chi_2+\chi_5$ if $p>3$). If $p>2$ there is a $G$-invariant complement $Q_{a'}$ for $Q_a$, giving a regular map ${\cal N}=M_{a'}\in E_p(\M)$ with $\dim({\cal N})=3$ and $K\cong Q/Q_{a'}$ affording $\pi_a-\chi_1$ ($=\chi_4$ if $p>3$).

Since the octahedron is $2$-face-colourable (i.e.~its dual, the cube, embeds a bipartite graph), there is a $2$-dimensional $G$-submodule $P_b>P_1$ spanned by the sums of the two monochrome subsets $\Psi\subset\Phi$, and a $5$-dimensional $G$-submodule $P_a+P_b<P^1$; these lead, via $G$-submodules $Q_b=P_b/P_1$ and $Q_a\oplus Q_b$ of $Q$, to coverings $\M_b$ and $\M_{a,b}$ in $E_p(\M)$ of dimensions $c=6$ and $3$ respectively. If $p>3$ then since $P_b$ affords the character $\chi_1+\chi_2$ of $G$, these two coverings afford $\chi_4+\chi_5$ and $\chi_5$ respectively. There is also a $7$-dimensional $G$-submodule $P_{b'}<P$ consisting of the elements $\sum_{\phi\in\Phi}a_{\phi}\phi$ such that $\sum_{\phi\in\Psi}a_{\phi}=\sum_{\phi\in\Psi'}a_{\phi}$ where $\Psi$ and $\Psi'$ are the monochrome subsets; it contains $P_1$, so its image $Q_{b'}$ in $Q$ is a $6$-dimensional $G$-submodule leadng to a cyclic covering $\M_{b'}$  affording $\chi_2$. Finally, $Q_{a'}\cap Q_{b'}$ is a $3$-dimensional $G$-submodule of $Q$ (the reduction mod~$(p)$ of the natural representation $\rho_n$ of $G$) giving a $4$-dimensional covering $\M_{a',b'}$ affording $\chi_2+\chi_4$. Together with $\M_0$ this accounts for all seven maps in $E_p(\M)$ when $p>3$.

The direct sum decomposition $Q=Q_b\oplus Q_a\oplus (Q_{a'}\cap Q_{b'})$, with direct factors affording the irreducible representations $\rho_2$, $\rho_4$ and $\rho_5$ corresponding to the characters $\chi_2$, $\chi_4$ and $\chi_5$, is also valid when $p=3$. (The character $\chi_3$, which reduces mod~$(3)$ to $\chi_1+\chi_2$, is not a summand of $\chi$; the remaining four characters in Table~\ref{chartS4} correspond to irreducible characters over $\F_3$.) Thus $E_3(\M)$ consists of the seven regular maps described above.

{\color{red}
}


\section{The dodecahedron}\label{Dod}

We now consider coverings of the dodecahedron ${\cal M}=\{5, 3\}$ branched over its face-centres (or, dually, of the icosahedron $\{3, 5\}$, branched over its vertices). In this case $G=\Delta(3, 2, 5)\cong A_5$, acting on the set $\Phi$ of twelve faces of $\cal M$, so the stabilisers of faces are the Sylow $5$-subgroups $H\cong C_5$ of $G$. The maps ${\cal N}\in E_p(\M)$ have type $\{5p, 3\}$ and genus $5p^c-6p^{c-1}+1$ where $c=\dim({\cal N})$. 

In $A_5$ there are five conjugacy classes: the identity, fifteen double transpositions, twenty $3$-cycles, and two classes of twelve $5$-cycles. Hence there are five irreducible characters. The character table is as shown in Table~\ref{chartA5}, with $\lambda,\mu=(1\pm\sqrt{5})/2$.
\begin{table}[htb]
\centering
\begin{tabular}{|c|c|c|c|c|c|}
\hline
&$1$&$(..)(..)$&$(...)$&$(.....)^+$&$(.....)^-$\\
\hline
$\chi_1$&1&1&1&1&1\\
$\chi_2$&3&-1&0&$\lambda$&$\mu$\\
$\chi_3$&3&-1&0&$\mu$&$\lambda$\\
$\chi_4$&4&0&1&-1&-1\\
$\chi_5$&5&1&-1&0&0\\ \hline
\end{tabular}
\caption{The character table of
$A_5$.}\label{chartA5}
\end{table}
 In addition to the principal character $\chi_1$, there are algebraically conjugate characters $\chi_2$ and $\chi_3$ obtained from the natural representation $\rho_n$ of $G$ as a rotation group, while the irreducible characters $\chi_4$ and $\chi_5$ are the non-principal summands of the permutation characters corresponding to the doubly transitive natural permutations representations of $G$ as $A_5$ and as $PSL_2(5)$.

The subgroup $H$ consists of the identity element, and two elements each from the two conjugacy classes of $5$-cycles. By averaging their values over $H$, one sees that the characters $\chi_i$ have multiplicities $1, 1, 1, 0$ and $1$ in the permutation character $\pi$ on $P$, that is,
\[\pi=\chi_1+\chi_2+\chi_3+\chi_5,\]
so
\[\chi=\chi_2+\chi_3+\chi_5,\]
and hence $H_1(S;{\mathbb C})$ is a direct sum of irreducible submodules of dimensions $3, 3$ and $5$.

The primes dividing $|G|$ are $2, 3$ and $5$, so we can use this decomposition to find the decomposition of $H_1(S;{\mathbb F}_p)$ for any prime $p>5$. If $p\equiv\pm 1$ mod~$(5)$, so that $5$ has a square root in ${\mathbb F}_p$, then $H_1(S;{\mathbb F}_p)$ is also a direct sum of irreducible submodules of dimensions $3, 3$ and $5$; in this case we obtain seven coverings ${\cal N}\in E_p(\M)$, of dimensions $c=3, 3, 5, 6, 8, 8$ and $11$. If $p\equiv\pm 2$ mod~$(5)$, on the other hand, $5$ has no square root in ${\mathbb F}_p$; in this case $H_1(S;{\mathbb F}_p)$ is a direct sum of irreducible submodules of dimensions $5$ and $6$, and there are three coverings ${\cal N}\in E_p(\M)$, with $c=5, 6$ and $11$. In either case, the first and second of these three maps are $\M_{a'}$ and $\M_a$, arising from the antipodal symmetry of the dodecahedron, and the third is $\M_0$; these coverings afford the characters $\chi_5$, $\chi_2+\chi_3$ and $\chi$ respectively. When $p\equiv\pm 1$ mod~$(5)$ there are two $3$-dimensional submodules $Q_n$ and $Q_{n*}$, obtained from the natural representation $\rho_n$ and its algebraic conjugate $\rho_{n*}$ of $G$; these lead to maps $\M_n$ and $\M_{n*}$ in $E_p(\M)$ with $c=8$ and $K$ affording $\chi_2+\chi_5$ and $\chi_3+\chi_5$. Finally, in this case there are also two $8$-dimensional submodules $Q_a\oplus Q_n$ and $Q_a\oplus Q_{n*}$, leading to $3$-dimensional coverings $\M_{a,n}$ and $\M_{a,n*}$ affording $\chi_3$ and $\chi_2$.

The direct sum decomposition $Q=Q_a\oplus Q_{a'}$, with direct factors affording the characters $\chi_5$ and $\chi_2+\chi_3$, is also valid when $p=3$ or $5$. Over $\F_3$, $\chi_5$ splits as $\chi_1+\chi_4$, while $\chi_2+\chi_3$ is irreducible (but splits into two algebraically conjugate characters over $\F_9$).


{\color{red}


}

\section{The icosahedron}\label{Ico}

We now consider coverings of the icosahedron ${\cal M}=\{3, 5\}$ branched over its face-centres (or, dually, of the dodecahedron $\{5, 3\}$, branched over its vertices). The maps ${\cal N}\in E_p(\M)$ have type $\{3p, 5\}$ and genus $9p^c-10p^{c-1}+1$, where $c=\dim({\cal N})$.

As in the case of the dodecahedron we have $G\cong A_5$, but now acting on the set $\Phi$ of twenty faces of $\cal M$, so the stabilisers of faces are the ten Sylow $3$-subgroups $H\cong C_3$ of $G$. The conjugacy classes and character table of $G$ are as before, but now $H$ consists of the identity and two $3$-cycles. In this case we find that
\[\pi=\chi_1+\chi_2+\chi_3+2\chi_4+\chi_5,\]
so
\[\chi=\chi_2+\chi_3+2\chi_4+\chi_5.\]
Thus $H_1(S;{\mathbb C})$ is a direct sum of irreducible submodules of dimensions $3, 3, 4$  (with multiplicity $2$) and $5$. By contrast with the previous examples, the existence of a summand with multiplicity greater than $1$ implies that the direct sum decompositions of $\C\Phi$ and of $H_1(S;\C)$ are not unique: in each case the $8$-dimensional submodule affording the character $2\chi_4$ contains infinitely many mutually isomorphic $4$-dimensional irreducible submodules affording $\chi_4$ (see \S 2.5).

As in the case of the dodecahedron, we can use this to obtain the decomposition of $H_1(S;{\mathbb F}_p)$ for primes $p>5$, but first we need to consider this lack of uniqueness of $4$-dimensional submodules.

There is a chiral pair of $G$-invariant equivalence relations $\sim_b$ and $\sim_c$ on $\Phi$, each with five classes of four faces, corresponding to the inclusion of $H$ in two subgroups of $G$ isomorphic to $A_4$. Their equivalence classes $\Psi$ provide basis elements $\underline\Psi$ for a pair of $5$-dimensional submodules $P_b, P_c>P_1$ of $P$, and hence a pair of $4$-dimensional submodules $Q_b$ and $Q_c$ of $Q$. There is an isomorphism $Q_b\to Q_c,\, q\mapsto q'$, induced by the antipodal automorphism $i$ of $\M$, and these submodules each afford the character $\chi_4$ which has multiplicity $2$ in the decomposition of $\chi$. They generate an $8$-dimensional submodule $Q_{b,c}=Q_b\oplus Q_c$ which contains $p-1$ other submodules isomorphic to them, each of the form $Q(\lambda)=\{q+\lambda q'\mid q\in Q_b\}$ for some $\lambda\in\F_p^*$. Now $i$ tranposes $\sim_b$ and $\sim_c$, and hence also transposes $Q(0):=Q_b$ and $Q(\infty):=Q_c$, so it permutes these $p+1$ submodules $Q(\lambda)$ ($\lambda\in{\mathbb P}^1(p)=\F_p\cup\{\infty\}$), leaving $Q(1)$ and $Q(-1)$ (if $p>2$) invariant and otherwise transposing pairs $Q(\lambda^{\pm 1})$. Since $A=G\times\langle i\rangle$, the $G$-submodules $Q(\pm 1)$ are actually $A$-submodules, affording the irreducible characters $\chi_4\otimes\varphi_{\pm 1}$ of $A$, where $\varphi_1$ and $\varphi_{-1}$ are the principal and non-principal characters of $\langle i\rangle\cong C_2$. 

If $p\equiv\pm 1$ mod~$(5)$ then $H_1(S;{\mathbb F}_p)$ is a direct sum of irreducible submodules of the same dimensions and multiplicities as those for $H_1(S;\C)$. The two irreducible submodules $Q_n$ and $Q_{n*}$ of dimension $3$ are obtained from the natural representation of $G$, those of dimension $4$ are the submodules $Q(\lambda)<Q_{a,b}$ defined above, and there is also one irreducible submodule $Q_5$ of dimension $5$.
There are eight submodules intersecting $Q_{b,c}$ trivially, namely the various sums of $Q_n$, $Q_{n*}$ and $Q_5$, with codimensions $8, 11, 11, 13, 14, 16, 16$ and $19$; by taking the direct sum of one of these with a submodule $Q(\lambda)$ or $Q_{a,b}$ we obtain $8(p+1)$ submodules with codimensions $4, 7, 7, 9, 10, 12, 12$ and $15$, or eight with codimensions $0, 3, 3, 5, 6, 8, 8$ and $11$. Excluding the submodule $Q$ from the last eight, this gives a total of $8p+23$ proper submodules, and hence the same number of coverings ${\cal N}\in E_p(\M)$, with $\dim({\cal N})=3, 4, \ldots, 16$ and $19$. Of these maps, $8(p-1)$ are chiral, namely those corresponding to submodules containing a single submodule $Q(\lambda)$, with $\lambda\ne\pm 1$, while the remaining $8+16+7=31$ are regular.

If $3<p\equiv\pm 2$ mod~$(5)$ then $H_1(S;{\mathbb F}_p)$ is the direct sum of four irreducible submodules: two are isomorphic $4$-dimensional submodules of the form $Q(\lambda)$, and the other two are $Q_5$ of dimension $5$, and $Q_{b,c}$ of dimension $6$. There are four submodules intersecting $Q_{a,b}$ trivially, namely $Q_{b,c}\oplus Q_5$, $Q_{b,c}$, $Q_5$ and $0$ with codimensions $8, 13, 14$ and $19$; by taking the direct sums of these with a submodule $Q(\lambda)$ or $Q_{a,b}$ we obtain $4(p+1)$ submodules with codimensions $4, 9, 10$ and $15$, or four with codimensions $0, 5, 6$ and $11$. Excluding $Q$ as before, we obtain $4p+11$ proper submodules, and the same number of coverings ${\cal N}\in E_p(\M)$, with $\dim({\cal N})=4, 5, 6, 8, 9, 10, 11, 13, 14, 15$ and $19$. These include $4(p-1)$ chiral maps, corresponding to submodules containing a single $Q(\lambda)$, with $\lambda\ne\pm 1$, while the other $15$ are regular.


For any $p$ the antipodal symmetry of the icosahedron gives a $9$-dimensional submodule $Q_a=Q(1)\oplus Q_5$ affording $\chi_4+\chi_5$, leading to a $10$-dimensional covering $\M_a$ affording $\chi_2+\chi_3+\chi_4$; if $p>2$ there is also a complementary $10$-dimensional submodule $Q_{a'}=Q_n\oplus Q_{n*}\oplus Q(-1)$ affording $\chi_2+\chi_3+\chi_4$, leading to a $9$-dimensional covering $\M_{a'}$ affording $\chi_4+\chi_5$.  The remaining summands $\chi_1, \chi_2, \chi_3$ and $\chi_5$ of $\pi|_G$ each extend to an irreducible character of $A$, so that $\pi|_A$ is a sum of six distinct irreducible characters, and thus $A$ has rank $6$ on $\Phi$.

{\color{red}


}

\section{The dihedron}\label{Dih}

If $\M$ is the dihedron $\{n, 2\}$, then
\[G=\Delta(2, 2, n)=\langle x, y, z\mid x^2=y^2=z^n=xyz=1\rangle.\]
This can be identified with the dihedral group
\[D_n=\langle a, b\mid a^n=b^2=1, a^b=a^{-1}\rangle,\]
where $x=b$ and $y=ab$ are rotations through $\pi$, fixing a vertex and edge-centre respectively and transposing the two faces, and $z=a$ is a rotation through $2\pi/n$, preserving the faces, so that $H=\langle a\rangle\cong C_n$. Then $H_1(S;\Z)$ is an infinite cyclic group, with $a$ acting on it as the identity and $b$ inverting every element. Thus
\[\pi = \chi_1+\chi_2\]
and hence
\[\chi=\chi_2,\]
where $\chi_1$ is the principal character of $G$ and $\chi_2(g)=(-1)^k$ for each $g=a^jb^k\in G$. It follows that for each integer $l\ge 1$ there is a single regular covering $\cal N$ of exponent $l$, with a cyclic covering group, corresponding to the unique subgroup of index $l$ in $H_1(S;\Z)$. This map is the dihedron $\{ln, 2\}$, a regular map of genus $0$ with $G\cong D_{ln}$ and $A=G\times\langle c\rangle\cong D_{ln}\times C_2$, where $c$ is the reflection fixing all the vertices and edges of $\cal N$ and transposing its two faces.

\section{The hosohedron}\label{Hos}

If $\M$ is the hosohedron $\{2, n\}$, then for each prime $p$ the maps ${\cal N}\in E_p(\M)$ are all regular, of type $\{2p,n\}$ and genus
\[1+\frac{p^{c-1}}{2}\bigl(n(p-1)-2p\bigr),\]
where $c=\dim({\cal N})$. As in the case of the dihedron, $G$ is isomorphic to $D_n$, but now acting naturally with degree $n$ on $\Phi$, so that the subgroup stabilising a face is now $H=\langle b\rangle\cong C_2$. This makes the analysis of the coverings much more complicated than before.

\subsection {Conjugacy classes and characters of $D_n$}

First suppose that $n$ is odd. Then apart from the conjugacy class $\{1\}$, there are $(n-1)/2$ conjugacy classes $\{a^{\pm j}\}$ for $j=1, 2, \ldots, (n-1)/2$, and a single conjugacy class $\{a^jb \mid j=0, 1, \ldots, n-1\}$ consisting of the involutions in $G$.

In addition to the characters $\chi_1$ and $\chi_2$ of degree $1$ defined in the preceding section, $G$ has $(n-1)/2$ irreducible characters $\xi_k$ of degree $2$, for $k=1, 2, \ldots, (n-1)/2$, where $\xi_k(a^{\pm j})=\alpha_{jk}:=\zeta_n^{jk}+\zeta_n^{-jk}$ with $\zeta_n=\exp(2\pi i/n)$. The character table of $G$ is as in Table~\ref{chartDnodd}, where the column headed $a^{\pm j}$ and the row labelled $\xi_k$ represent  $(n-1)/2$ conjugacy classes and characters respectively, as $j$ and $k$ each range over $\{1, 2, \ldots, (n-1)/2\}$. For each $k$, the kernel of the representation $\rho_k$ corresponding to $\xi_k$ is the unique subgroup $\langle a^m\rangle$ of order $\gcd(k, n)$ in $G_0:=\langle a\rangle$, where $m=n/\gcd(k, n)\ge 3$, so $\rho_k$ can be regarded as a faithful representation of $G/\langle a^m\rangle\cong D_m$.
 
\begin{table}[htb]
\centering
\begin{tabular}{|c|c|c|c|}
\hline
&$1$&$a^{\pm j}$&$a^jb$\\
\hline
$\chi_1$&1&1&1\\
$\chi_2$&1&1&$-1$\\
$\xi_k$&2&$\alpha_{jk}$&0\\ \hline
\end{tabular}
\caption{The character table of $D_n$ for odd $n$.}\label{chartDnodd}
\end{table}

By computing the inner product of each irreducible character with $\pi$ we find that
\[\pi=\chi_1+\sum_{k=1}^{(n-1)/2}\xi_k\]
and hence
\[\chi=\sum_{k=1}^{(n-1)/2}\xi_k.\]

If $n$ is even the conjugacy classes of $G$ are $\{1\}$, $\{a^{n/2}\}$, $\{a^{\pm j}\}$ for $j=1, 2, \ldots, (n-2)/2$, $\{a^jb\mid j\;\hbox{is even}\}$ and $\{a^jb\mid j\;\hbox{is odd}\}$. In addition to $\chi_1$ and $\chi_2$, defined as above for $n$ odd, there are two more characters $\chi_3$ and $\chi_4$ of degree $1$, together with irreducible characters $\xi_k$ of degree $2$ for $k=1, 2, \ldots, (n-2)/2$. The character table is as in Table~\ref{chartDneven}, with the columns headed $a^{\pm j}$ and $a^jb$ representing $1+(n-2)/2=n/2$ and $2$ conjugacy classes respectively, and the row labelled $\xi_k$ representing $(n-2)/2$ characters. We can take $H=\langle b\rangle$, so that
\[\pi=\chi_1+\chi_3+\sum_{k=1}^{(n-2)/2}\xi_k\]
and hence
\[\chi=\chi_3+\sum_{k=1}^{(n-2)/2}\xi_k.\]

 \begin{table}[htb]
\centering
\begin{tabular}{|c|c|c|c|}
\hline
&$1$&$a^{\pm j}$&$a^jb$\\
\hline
$\chi_1$&1&1&1\\
$\chi_2$&1&1&$-1$\\
$\chi_3$&1&$(-1)^j$&$(-1)^j$\\
$\chi_4$&1&$(-1)^j$&$(-1)^{j+1}$\\
$\xi_k$&2&$\alpha_{jk}$&0\\ \hline
\end{tabular}
\caption{The character table of $D_n$ for even $n$.}\label{chartDneven}
\end{table}

\subsection{Submodule structure}

Instead of considering the reduction mod~$(p)$ of the characters appearing in $\pi$, as in the preceding cases, it is more convenient to use the fact that the group $G=D_n$ has a subgroup $G_0=\langle a\rangle$ acting regularly on $\Phi$. This allows us to give a more explicit description of the $G$-submodules of $P$ and $Q$, and hence of the maps ${\cal N}\in E_p(\M)$. We will assume that $p$ is coprime to $n$.

In order to find the decompositions of $P=\F_p\Phi$, and hence of $H_1(S;\F_p)\cong Q=P/P_1$, as $G$-modules, we first consider them as $G_0$-modules. Since $G_0$ acts regularly on $\Phi$, one can identify $P$ with the group algebra $\F_pG_0$ of $G_0$ over $\F_p$, or equivalently with the polynomial algebra $\F_p[x]/(x^n-1)$, the elements $a^j\in G_0$ corresponding to the images of the powers $x^j$ of $x$. The $G_0$-submodules are then the ideals of this algebra, each generated by the image of the ideal in $\F_p[x]$ generated by some polynomial dividing $x^n-1$. To find these polynomials we need to determine the irreducible factors of $x^n-1$ in $\F_p[x]$.

In $\Z[x]$ we have a factorisation
\begin{equation}
x^n-1=\prod_{m|n}\Phi_m(x),
\end{equation}
where $\Phi_m(x)$ in the $m$th cyclotomic polynomial, the minimal polynomial of the primitive $m$th roots of $1$ in $\C$. In fact, equation~(1), together with the equation $\Phi_1(x)=x-1$, can be regarded as a recursive definition of $\Phi_n(x)$ in terms of the polynomials $\Phi_m(x)$ for the proper divisors $m$ of $n$.

Each cyclotomic polynomial is irreducible in $\Z[x]$ (in fact, in $\Q[x]$), but in almost all cases it factorises when reduced mod~$(p)$, that is, when regarded as a polynomial in $\F_p[x]$. We are assuming that $n$ is coprime to $p$, so the set $\Omega=\Omega_n$ of $n$th roots of $1$ in $\overline\F_p$ contains $n$ elements. This set is partitioned into disjoint subsets $\Pi_m$ ($m|n$) consisting of the $\varphi(m)$ primitive $m$th roots of $1$ in $\overline\F_p$. For each divisor $m$ of $n$, let $e=e_m$ denote the multiplicative order of $p$ as a unit mod~$(m)$, so that $m$ divides $p^e-1$ but not $p^f-1$ for any $f<e$. Then each $\omega\in\Pi_m$ generates the subfield $\F_{p^e}$ of $\overline\F_p$. The restriction to $\F_{p^e}$ of the Frobenius automorphism $\phi: t\mapsto t^p$ of $\overline\F_p$ generates the Galois group $C={\rm Gal}\,\F_{p^e}\cong C_e$ of this subfield, which partitions $\Pi_m$ into $\varphi(m)/e$ orbits $\Gamma=\{\omega, \omega^p, \ldots, \omega^{p^{e-1}}\}$ of length $e$. Each such orbit $\Gamma$ is the set of roots $\omega^{p^i}$ of an irreducible monic factor $f^{\Gamma}(x)$ of degree $e$ of $\Phi_m(x)$ in $\F_p[x]$. As $\Gamma$ ranges over all the orbits of $C$ on $\Omega$, these polynomials $f^{\Gamma}(x)$ give all the irreducible factors of $x^n-1$, and the ideals $P^{\Gamma}$ they generate give all the maximal $G_0$-submodules of $P$. By taking the ideals generated by all products of these irreducible factors we obtain all the $G_0$-submodules of $P$ as intersections of these maximal submodules. In particular, $P$ is the direct sum of the irreducible $G_0$-submodules $P_{\Gamma}\cong P/P^{\Gamma}$ generated by the complementary factors $f_{\Gamma}(x)=(x^n-1)/f^{\Gamma}(x)$ of $x^n-1$, with the generator $a$ of $G_0$ having minimal polynomial $f^{\Gamma}(x)$ on $P_{\Gamma}$. For instance, the orbit $\Gamma=\Pi_1=\{1\}$ corresponds to the submodule $P_1$ fixed by $G_0$; this is generated by $(x^n-1)/(x-1)=x^{n-1}+\cdots+x+1$, the product of all the irreducible factors except $\Phi_1(x)=x-1$.  More generally, if $m$ divides $\gcd(n, p-1)$ then $\Pi_m$ consists of $\varphi(m)$ singleton orbits $\Gamma=\{\omega\}$, each giving a linear factor $f^{\Gamma}(x)=x-\omega$ of $\Phi_n(x)$ and corresponding to a $1$-dimensional direct summand $P_{\Gamma}$ of $P$ on which $a$ has eigenvalue $\omega$.


To understand the structures of $P$ and $Q$ as $G$-modules, we need to consider the action of $b$, which acts on $G_0$ by conjugation as inversion $\iota: a^j\mapsto a^{-j}$. Given any orbit $\Gamma\subseteq \Pi_m$ of the Galois group $C=\langle\phi\rangle$ on $\Omega$, the inverses of its elements also form such an orbit $\Gamma^*$, so that $\Delta=\Gamma\cup\Gamma^*$ is an orbit of the group $B=\langle\phi,\iota\rangle$. The orbit $\Gamma$ containing $\omega$ coincides with $\Gamma^*$ if and only if $\omega^{p^f}=\omega^{-1}$, or equivalently, $p^f\equiv -1$ mod~$(m)$, for some integer $f$. (This is always true if $m=1$ or $2$, but otherwise it implies that $e$ is even and $f\equiv e/2$ mod~$(e)$.) This condition depends only on $p$ and $m$, so for a given $m$, either all of the orbits $\Gamma\subseteq\Pi_m$ of $C$ are self-inverse, each forming an orbit $\Delta$ of $B$, or none of them are, so that each orbit $\Delta$ of $B$ is a union of two mutually inverse orbits $\Gamma$ and $\Gamma^*$ of $C$. In the first case, the $G_0$-submodules $P^{\Gamma}$ and $P_{\Gamma}$ of $P$ are also $G$-modules, with $P_{\Gamma}\cong P/P^{\Gamma}$ irreducible and $e$-dimensional, whereas in the second case $P_{\Gamma}\oplus P_{\Gamma^*}$ is an irreducible $2e$-dimensional $G$-submodule, isomorphic to $P/(P^{\Gamma}\cap P^{\Gamma^*})$, with $b$ transposing the two direct factors. For instance, any $1$-dimensional $G_0$-submodules $P_{\Gamma}$ are of the second type, apart those with $m=1$ or $2$. In either case, we will denote this irreducible submodule by $R_{\Delta}$, and we let $d=e$ or $2e$ denote its dimension. The kernel of the representation of $G$ on $R_{\Delta}$ is the normal subgroup $\langle a^m\rangle\cong C_{n/m}$, so that $G$ acts on $R_{\Delta}$ as $G/\langle a^m\rangle\cong D_m$.

These $G$-submodules $R_{\Delta}$ are the irreducible summands in the direct sum decomposition of the $G$-module $P=P_1\oplus P^1$; deleting the $1$-dimensional submodule $P_1$ gives the corresponding decomposition of $Q=P/P_1\cong P^1$. We therefore know all the proper $G$-submodules $L$ of $Q$ and hence all the maps $\cal N$ in $E_p(\M)$. The possible values of $c=\dim({\cal N})$ range from $1$ (corresponding to $L=P^{\Delta}/P_1$ where $n$ is even and $\Delta=\Pi_2=\{-1\}$) to $n-1$ (for ${\cal N}=\M_0$, corresponding to $L=0$). The number $\nu$ of irreducible summands $R_{\Delta}$ in the direct sum decomposition of the $G$-module $Q$ is equal to the number of orbits $\Delta\ne\{1\}$ of $B$ on $\Omega$. Since these summands are mutually non-isomorphic, the total number of $G$-submodules of $Q$ is $2^{\nu}$; thus $Q$ has $2^{\nu}-1$ proper $G$-submodules and hence $|E_p(\M)|=2^{\nu}-1$.

\subsection{Examples}

To illustrate the preceding arguments, suppose that $n=95$ and $p=7$.  We consider the divisors $m=1, 5, 19$ and $95$ of $n$ in turn.

The divisor $m=1$ leads to a single orbit $\Delta=\Pi_1=\{1\}$ of $B$, with the corresponding $1$-dimensional direct summand $R_{\Delta}=P_1$ of $P$, generated by $(x^{95}-1)/(x-1)$ and affording the principal representation of $G$.

Next we consider $m=5$. The prime $p=7$ has order $e=4$ in ${\mathbb Z}^*_5$; since $\varphi(5)=4$, the Galois group $C={\rm Gal}\,\F_{7^4}\cong C_4$ has $\varphi(5)/e=1$ orbit $\Gamma$ on the set $\Pi_5\subset\F_{7^4}\subset\overline{\mathbb F}_7$ of primitive $5$th roots of $1$. The cyclotomic polynomial $\Phi_5(x)=f^{\Gamma}(x)=x^4+x^3+x^2+x+1$ is therefore irreducible in $\F_7[x]$, giving a $4$-dimensional irreducible $G_0$-submodule $P_{\Gamma}=(x^{95}-1)/(f^{\Gamma}(x))$ of $P$. Since $-1\in\langle 7\rangle$ in ${\mathbb Z}^*_5$, this is also a $4$-dimensional irreducible $G$-submodule $R_{\Delta}$, with $G$ acting on it as $D_5$. When we factor out $P_1$, this appears as a direct summand of $Q=H_1(S;\Z_7)$.

Now let $m=19$. Since $7$ has order $e=3$ in ${\mathbb Z}^*_{19}$, and there are $\varphi(19)=18$ primitive $19$th roots of $1$ in $\F_{7^3}$, forming six orbits $\Gamma$ under the Galois group $C\cong C_3$ of this field, we obtain six $3$-dimensional irreducible $G_0$-submodules $P_{\Gamma}$ as summands of $P$, corresponding to six irreducible cubic factors $f^{\Gamma}(x)$ of $\Phi_{19}(x)$ in $\F_7[x]$. Since $-1\not\in\langle 7\rangle$ in ${\mathbb Z}^*_{19}$, these six orbits form three mutually inverse pairs $\Gamma$ and $\Gamma^*$, so they merge to give three orbits $\Delta=\Gamma\cup\Gamma^*$ of $B$. These correspond to three $6$-dimensional irreducible $G$-submodules $R_{\Delta}=P_{\Gamma}\oplus P_{\Gamma^*}$ of $P$, and hence also of $Q$, with $G$ acting on each as $D_{19}$.

Finally let $m=95$. Since $\Z_{95}\cong \Z_5\oplus\Z_{19}$ we see that $7$ has order $e={\rm lcm}\{4, 3\}=12$ in ${\mathbb Z}^*_{95}$, so $\Phi_{95}(x)$ is a product of $\varphi(95)/12=6$ irreducible polynomials $f^{\Gamma}(x)$ of degree $12$ in $\F_7[x]$, giving six $12$-dimensional  irreducible $G_0$-submodules $P_{\Gamma}$ as summands of $P$. Since $-1\not\in\langle 7\rangle$ in ${\mathbb Z}^*_{19}$ we have  $-1\not\in\langle 7\rangle$ in ${\mathbb Z}^*_{95}$, so these give three $24$-dimensional irreducible $G$-submodules $R_{\Delta}=P_{\Gamma}\oplus P_{\Gamma^*}$ of $P$, and hence of $Q$, with $G$ acting faithfully on each. 

Adding up the dimensions of all these irreducible submodules $R_{\Delta}$, we have
\[1+4+3\times 6+3\times 24=95=\dim P,\]
so they give the complete direct sum decomposition $\oplus_{\Delta}R_{\Delta}$ of $P$ as a $G$-module; deleting the $1$-dimensional submodule $P_1$ gives the corresponding decomposition of $Q$. Since $Q$ has $\nu=7$ mutually non-isomorphic irreducible summands, we obtain $2^7-1=127$ proper $G$-submodules $L<Q$. Thus $E_7(\M)$ consists of $127$ regular maps $\cal N$, with $\dim({\cal N})=4$ (one map), $6$ (three maps), and so on, up to $94$ (one map, $\M_0$).

Examples with even values of $n$ can be treated similarly, except that now there is always at least one $1$-dimensional summand $R_{\Delta}$, corresponding to the orbit $\Delta=\Gamma=\Pi_2=\{-1\}$ of $B$ and the polynomial $f^{\Gamma}(x)=\Phi_2(x)=x+1$. This corresponds to a cyclic covering $\cal N$ on which $a$ has eigenvalue $-1$ (see~\cite[\S 9(d)]{JSu}).

For instance, if $n=4$ then $E_p(\M)$ consists of regular maps of type $\{2p, 4\}$ and genus $g=1+p^{e-1}(p-2)$. If $p\ne 2$ then $B$ has orbits $\Delta=\{1\}, \{-1\}$ and $\{\zeta_4^{\pm 1}\}$ on $\Omega_4$, the last two providing irreducible $G$-submodules of $Q$ of dimensions $1$ and $2$, affording the characters $\chi_3$ and $\xi_1$ of $G$. Thus $E_p(\M)$ consists of three maps ${\cal N}$ with $c=\dim({\cal N})=1, 2$ and $3$, of genus $g=p-1, (p-1)^2$ and $(p-1)(p^2-p-1)$, affording the characters $\chi_3$, $\xi_1$ and $\chi_3+\xi_1$. The maps with $c=1$ lie on the Accola-Maclachlan surfaces with $8(g+1)$ automorphisms~\cite{Acc, Macl}; for $p=3,5,7,11$ they are the duals of R2.2, R4.4, R6.4 and R10.11 in~\cite{C}. Those with $c=2$ (the maps $\{2p,4\}_4$ in~\cite[\S8.6 and Table 8]{CM}) are the duals of R4.3, R16.3, R36.4 and R100.5 in~\cite{C} for these primes.


\subsection{Prime $n$}

If $n$ is prime then the only divisors $m$ of $n$ are $m=1$, giving rise to the $G$-submodule $P_1\le P$, and $m=n$. If $p\ne n$ then $\Phi_n(x)$ factorises in $\F_p[x]$ as a product of $(n-1)/e$ irreducible polynomials $f^{\Gamma}(x)$ of degree $e$ equal to the order of $p$ in the group $\Z_n^*\cong C_{n-1}$. As a $G_0$-module, $Q$ is a direct sum of $(n-1)/e$ corresponding irreducible submodules $Q_{\Gamma}$ of dimension $e$. If $n=2$ or if $e$ is even then $-1\in\langle p\rangle$ in $\Z_n^*$, so these are the $\nu=(n-1)/e$ irreducible summands $R_{\Delta}$ of dimension $d=e$, for $\Delta=\Gamma$, in the direct sum decomposition of $Q$ as a $G$-module. If $n>2$ and $e$ is odd then $-1\not\in\langle p\rangle$ in $\Z_n^*$, so we obtain $\nu=(n-1)/2e$ direct summands $R_{\Delta}=P_{\Gamma}\oplus P_{\Gamma^*}$ of dimension $d=2e$ in this decomposition. In either case, there are $2^{\nu}-1$ maps ${\cal N}\in E_p(\M)$, specifically ${\nu\choose i}$ with $c=\dim({\cal N})=id$ for each $i=1, 2, \ldots, \nu$, affording the sum of $i$ distinct irreducible characters of degree $d$. In particular, if $n>2$ and $p\equiv \pm 1$ mod~$(n)$ then $d=2$, so $\nu=(n-1)/2$, with the irreducible coverings affording the characters $\xi_k$ of $G$ for $k=1, 2, \ldots, (n-1)/2$.

(If $p=n$ then $Q\cong \F_p[x]/((x-1)^{n-1})=\F_n[x]/((x-1)^{n-1})$, with proper $G$-submodules of codimensions $c=1, 2,\ldots, n-1$, each generated by the image of $(x-1)^c$. It follows that $E_p(\M)$ consist of $n-1$ coverings, one for each of these dimensions $c$.)

In the next two subsections we give some simple illustrative examples.

\subsubsection{Example: $n=3$}

The first value of $n$ not covered by the section on the dihedron is $n=3$. In this case the maps ${\cal N}\in E_p(\M)$ are regular, of type $\{2p, 3\}$ and genus $1+\frac{1}{2}p^{c-1}(p-3)$ where $c=\dim({\cal N})$. If $p\ne 3$ then $p$ has order $e=1$ or $2$ in $\Z_3^*$ as $p\equiv 1$ or $-1$ mod~$(3)$, so as a $G_0$-module $Q$ is respectively a direct sum of two $1$-dimensional submodules, or is $2$-dimensional and irreducible. In either case $B$ has a single orbit $\Delta\ne\{1\}$ on $\Omega_3$, namely $\{\zeta_3^{\pm 1}\}$ where $\zeta_3$ is a primitive cube root of $1$, so $Q$ is an irreducible $2$-dimensional $G$-module; hence $E_p(\M)$ consists of a single map ${\cal N}=\M_0$ of genus $(p-1)(p-2)/2$, affording the character $\xi_1$ of $G$ (this is the map $\{2p,3\}_6$ in~\cite[\S8.6 and Table 8]{CM}). For instance, if $p=2$ this map $\cal N$ is the cube $\{4, 3\}$, with $\M={\cal N}/V_4$, and for $p=5, 7, 11$ or $13$ it is the dual of the regular map R6.1, R15.2, R45.3 or R66.1 in~\cite{C}.


\subsubsection{Example: $n=13$}

If $n=13$ the maps in ${\cal N}\in E_p(\M)$ have type $\{2p, 13\}$ and genus $1+\frac{1}{2}p^{c-1}(11p-13)$. If $p=\pm 1$ mod~$(13)$ then $e=1$ or $2$; in either case $d=2$, so there are $\nu=12/2=6$ irreducible coverings, affording the characters $\xi_1,\ldots, \xi_6$, and a total of $2^6-1=63$ coverings ${\cal N}\in E_p(\M)$, with $\dim({\cal N})=2, 4, \ldots, 12$. If $p\equiv 3$ or $9$ mod~$(13)$ then $e=3$; this is odd, so $d=6$, and hence there are $\nu=2$ irreducible coverings of dimension $6$; thus $E_p(\M)$ consists of these two irreducible coverings and $\M_0$ of dimension $12$. The same applies if $p\equiv\pm 4$ mod~$(13)$, so that $d=e=6$. If $p\equiv \pm 5$ mod~$(13)$ then $e=4$, so $d=4$ also; thus $\nu=3$, so $E_p(\M)$ consists of seven coverings of dimensions $4, 8$ and $12$. If $p\equiv\pm 2$ or $\pm 6$ mod~$(13)$ then $d=e=12$, giving one irreducible covering $\M_0$ of dimension $12$. 

\section{Branching over faces and vertices}

It is easy now to find the elementary coverings of the Platonic maps, where the branching is over their vertices rather than the faces: one simply uses the preceding results to find the coverings of the dual map $\M^*$ of $\M$, branched over its faces, and then takes the duals of the resulting maps.

One can also use similar methods to determine the elementary abelian coverings, where the branching is over both the faces and the vertices. In this case, we replace $S$ with the sphere $S'$ punctured at the set $\Phi'=\Phi\cup\Phi^*$, where $\Phi$ is the set of face-centres of $\M$ and $\Phi^*$ is its set of vertices (i.e.~the face-centres of $\M^*$). The resulting homology module $H_1(S';\F_p)$ for $G$ on $S'$ is the quotient $Q'=P'/P'_1$, where $P'$ is the permutation module $\F_p\Phi'=\F_p\Phi\oplus\F_p\Phi^*$ for $G$ on $\Phi'$, and $P_1'$ is the $1$-dimensional $G$-submodule of $P'$ spanned by $\underline\Phi'=\underline\Phi+\underline\Phi^*$. The representation and character of $G$ on $P'$ are $\rho'=\rho\oplus\rho^*$ and $\pi'=\pi+\pi^*$, where $\rho$ and $\rho^*$ are the representations on the permutation modules $P=\F_p\Phi$ and $P^*=\F_p\Phi^*$, and $\pi$ and $\pi^*$ are their associated characters. Since $G$ has two orbits on $\Phi'$, the principal character $\chi_1$ of $G$ has multiplicity $2$ in $\pi'$, afforded by the $1$-dimensional fixed submodules $P_1$ and $P_1^*$ of $P$ and $P^*$, spanned by $\underline\Phi$ and $\underline\Phi^*$; deleting one copy of $\chi_1$ gives the character
\[\chi'=\pi+\pi^*-\chi_1=\chi_1+\chi+\chi^*\]
of $G$ on $Q'$, where $\chi=\pi-\chi_1$ and $\chi^*=\pi^*-\chi_1$ are the characters on $Q=P/P_1$ and $Q^*=P^*/P^*_1$. We have already determined the $G$-module structures of $P$ and of $P^*$ (the latter in the context of $\M^*$ rather than $\M$), so we can immediately deduce the structure of $P'=P\oplus P^*$, using Lemma~2.1, and hence that of its quotient $Q'$. This gives us information about the set $E'_p(\M)$ of coverings of $\M$ by elementary $p$-groups, branched over faces and vertices, since these correspond to the proper $G$-submodules of $Q'$.

In all cases, $Q'$ is an extension of a fixed $G$-submodule $Q'_1=(P_1\oplus P^*_1)/P_1'\cong P_1$ ($\cong P_1^*$) by
\[P/(P_1\oplus P^*_1)
=(P\oplus P^*_1)/(P_1\oplus P^*_1)\oplus(P_1\oplus P^*)/(P_1\oplus P^*_1)
\cong Q\oplus Q^*.\]
If the modules $P$ and $P^*$ both split over their fixed submodules $P_1$ and $P_1^*$ (equivalently, if both $f=|\Phi|$ and $f^*=|\Phi^*|$ are coprime to $p$), then $Q'$ also splits over $Q'_1$: we have $P=P_1\oplus P^1$ and $P^*=P_1^*\oplus {P^*}^1$, so
\[P'=P_1\oplus P_1^*\oplus P^1\oplus {P^*}^1;\]
since $P'_1$ is a submodule of $P_1\oplus P^*_1$ there is a complement
\[Q'^1=(P'_1\oplus P^1\oplus {P^*}^1)/P'_1\cong P^1\oplus {P^*}^1\cong Q\oplus Q^*\]
for $Q'_1$ on $Q'$, giving
\[Q'=Q'_1\oplus Q'^1 \cong P_1\oplus Q\oplus Q^*.\]

Instead of dealing with all the cases in detail, we will consider one simple example, and just summarise the results in the remaining cases.

\subsection{The tetrahedron}

If $\M$ is the tetrahedron then $\M^*\cong\M$, and so $P\cong P^*$. We saw in \S3 that if $p\ne 2$ then $P=P_1\oplus P^1$, affording the representation $\rho_1\oplus\rho_4$ of $G$ with character $\chi_1+\chi_4$. It follows that
\[Q'=Q'_1\oplus Q_3\oplus Q_3^*,\]
affording the representation $\rho_1\oplus\rho_4\oplus\rho_4$ with character
\[\chi'=\chi_1+2\chi_4;\]
here $Q'_1$ affords $\rho_1$, while $Q_3$ and $Q_3^*$, the images in $Q'$ of $P^1$ and ${P^*}^1$, both afford $\rho_4$. Now $Q_3\cong Q_3^*\cong Q$, so $Q_3\oplus Q_3^*$ contains $p+1$ submodules $Q(\lambda)\cong Q$, where $\lambda\in{\mathbb P}^1(p)$ (see \S 2.5, and also \S7 for a similar phenomenon concerning the icosahedron). It follows that $Q'$ has the following proper submodules $L$: one submodule $Q_3\oplus Q_3^*$ of codimension $c=1$, $p+1$ submodules $Q'_1\oplus Q(\lambda)$ of codimension $3$, $p+1$ submodules $Q(\lambda)$ of codimension $4$, one submodule $Q'_1$ of codimension $6$, and one submodule $0$ of codimension $7$. Thus there are $|E_p'(\M)|=2p+5$ coverings $\cal N$ of $\M$, of these dimensions $c$, affording the representations $\rho_1$, $\rho_4$, $\rho_1\oplus \rho_4$, $\rho_4\oplus\rho_4$, and $\rho_1\oplus\rho_4\oplus\rho_4$ respectively. As in \S3, the maps ${\cal N}\in E_p'(\M)$ are all regular.

The coverings of dimension $3$ corresponding to $L=Q'_1\oplus Q^*_3$ and $Q'_1\oplus Q_3$ are respectively  branched over only the faces and only the vertices; the first of these is the covering $\M_0$ of type $\{3p, 3\}$ and genus $p^3-2p^2+1$ described in \S3, and the second is its dual. The remaining $2p+3$ coverings have branching of order $p^c-p^{c-1}$ at the four vertices and the four faces of $\M$, so they have type $\{3p, 3p\}$ and genus $3p^c-4p^{c-1}+1$. When $L=Q_3\oplus Q^*_3$ these are cyclic coverings of genus $3(p-1)$, obtained by assigning mutually inverse monodromy permutations in $K=C_p$ to the vertices and face-centres of $\M$; for instance, if $p=3, 5, 7, 11$ these are the maps R6.9, R12.8, R18.10, R30.10 in~\cite{C}. The $p-1$ coverings with $L=Q'_1\oplus Q(\lambda)$ and $\lambda\ne 0, \infty$ have genus $3p^3-4p^2+1$; when $p=3$ these are the maps R46.27 of type $\{9, 9\}_4$ and R46.28 of type $\{9, 9\}_{12}$ in~\cite{C}.

{\color{red}
}

\subsection{The remaining cases}

Here we simply give the homology character $\chi'$ for each of the remaining Platonic maps, in terms of the irreducible characters in the tables given earlier, together with a few comments. The $G$-module structure of $Q'$ and hence the properties of the elementary abelian coverings $\cal N$ of $\M$ can be deduced by methods similar to those used for the tetrahedron.

If $\M$ is the cube or octahedron then
\[\chi'=\chi_1+\chi_2+\chi_3+\chi_4+2\chi_5.\]
The presence of the $1$-dimensional rational characters $\chi_1$ and $\chi_2$ leads to both central and non-central cyclic coverings $\cal N$ of $\M$. As in the case of the tetrahedron, the coefficient $2$ of $\chi_5$ means that the number of coverings is unbounded as a function of $p$: in fact, since $\rho_5\oplus\rho_5$ has $p+1$ summands isomorphic to $\rho_5$ we have $|E_p'(\M)|=16p+47$ for $p>3$.



If $\M$ is the dodecahedron or icosahedron then
\[\chi'=\chi_1+2\chi_2+2\chi_3+2\chi_4+2\chi_5.\]
The character $\chi_1$ again leads to central cyclic coverings, but now there are no non-central cyclic coverings. The characters $\chi_4$ and $\chi_5$ are both rational, so the representations $\rho_4$ and $\rho_5$ are realised over $\F_p$. The same applies to $\rho_2$ and $\rho_3$ if $p\equiv\pm 1$ mod~$(5)$, whereas if $p\equiv\pm 2$ mod~$(5)$ then we obtain a $6$-dimensional irreducible representation over $\F_p$ with character $\chi_2+\chi_3$. The presence of four characters with coefficient $2$ leads to the conclusion that $|E_p'(\M)|\sim 2p^4$ as $p\to\infty$, the dominating terms arising from $G$-submodules of $Q'$ which are direct sums of four or five non-isomorphic irreducible submodules.


If $\M$ is the dihedron $\{n, 2\}$ or the hosohedron $\{2, n\}$ then
\[\chi'=\chi_1+\chi_2+\sum_{k=1}^{(n-1)/2}\xi_k\]
if $n$ is odd, and
\[\chi'=\chi_1+\chi_2+\chi_3+\sum_{k=1}^{(n-2)/2}\xi_k\]
if $n$ is even. The $1$-dimensional rational characters $\chi_k$ give rise to cyclic coverings, central when $k=1$ and non-central when $k>1$. The fact that $\chi'$ is multiplicity-free means that $|E_p'(\M)|$ is bounded above as a function of $p$, by $2^{(n+3)/2}-1$ or $2^{(n+4)/2}-1$ as $n$ is odd or even. These bounds are attained when $p\equiv 1$ mod~$(n)$, so that the representations corresponding to the characters $\xi_k$ all reduce mod~$(p)$ to $2$-dimensional representations over $\F_p$. For each $n$ there are infinitely many such primes $p$, by Dirichlet's theorem on primes in arithmetic progressions.

\section{Branching over edges}

Hypermaps (Grothendieck's {\em dessins d'enfants}, see~\cite{Gro, JSi2}) are generalisations of maps, which correspond to triangle groups $\Delta(p, q, r)$ with $q$ not necessarily equal to $2$. These arise if we extend the investigation of coverings of the Platonic maps to allow branching over their edges. The general principles are the same as for branching over the faces, except that $\Phi$ is now replaced with the set $\Phi^{\dagger}$ of edges of $\M$, the sphere is punctured at their midpoints, and $H$ is the subgroup $\langle y\rangle\cong C_2$ of $G$ leaving invariant an edge. Dual pairs of maps $\cal M$ can be treated together, since duality induces $G$-isomorphisms between their edge sets.

The decompositions of the resulting character $\chi^{\dagger}$ of $G$ on the homology module are as follows, with the notation for irreducible characters as before. For the tetrahedron (see \S 11.1 for applications),
\[\chi^{\dagger}=\chi_2+\chi_3+\chi_4.\]
For the cube and the octahedron,
\[\chi^{\dagger}=\chi_3+2\chi_4+\chi_5.\]
For the dodecahedron and the icosahedron,
\[\chi^{\dagger}=\chi_2+\chi_3+2\chi_4+3\chi_5.\]
For the dihedron and the hosohedron, if $n$ is odd then
\[\chi^{\dagger}=\sum_{k=1}^{(n-1)/2}\xi_k,\]
the same decomposition as for branching over the faces of the hosohedron, by the duality between edges and faces of this map; if $n$ is even then taking $H=\langle ab\rangle$ as the subgroup stabilising an edge (since we took $\langle b\rangle$ in \S 9.1 as the stabiliser of a face), we have
\[\chi^{\dagger}=\chi_1+\chi_4+\sum_{k=1}^{(n-2)/2} \xi_k.\]

\subsection{The tetrahedron}

Let $\M$ be the tetrahedron, so that the coverings $\cal N$ are hypermaps of type $(3, 2p, 3)$ and genus $2p^c-3p^{c-1}+1$. If $p\equiv 1$ mod~$(3)$ then $\chi^{\dagger}$ decomposes over $\F_p$ as above, with irreducible summands $\chi_2$, $\chi_3$ and $\chi_4$ of degrees $1$, $1$ and $3$, so we obtain seven coverings $\cal N$, with $c=1, 1, 2, 3, 4, 4$ and $5$. Since $\chi_2$ and $\chi_3$ are transposed by orientation-reversing automorphisms of $\M$, those with $c=1$ or $4$ occur in chiral pairs, whereas the other three are regular. For instance, if $p=7$ then the two cyclic coverings (with $c=1$) are duals of the chiral pair of hypermaps CH12.1 of type $(3,3,14)$ in~\cite{C}. If $2<p\equiv 2$ mod~$(3)$ then $\chi^{\dagger}$ decomposes as a sum of two irreducible characters $\chi_2+\chi_3$ and $\chi_4$ of degrees $2$ and $3$, so in this case there are three coverings $\cal N$ with $c=2, 3$ and $5$, all of them regular. For instance, if $p=5$ then the covering with $c=2$ is a dual of the hypermap RPH36.3 in~\cite{C}.


\section{Branching over vertices, edges and faces}

Just as in \S 10 we considered branching over the vertices and faces, one can also consider coverings which are branched over any combination of vertices, edges and faces. In the most general case, where we allow branching over all three sets, $G$ has three orbits on the punctures of the sphere; the permutation module $P''$ is therefore the direct sum of $P$, $P^*$ and a third permutation module $P^{\dagger}$ with $\Phi^{\dagger}$ as its basis, and the corresponding homology character is
\[\chi''=\pi+\pi^*+\pi^{\dagger}-\chi_1=2\chi_1+\chi+\chi^*+\chi^{\dagger}.\]
If $\M$ is the tetrahedron, for example, then $\chi''=2\chi_1+\chi_2+\chi_3+3\chi_4$.

The corresponding coverings of each $\cal M$ can now be found by using the techniques described earlier. For instance, in all cases the character $\chi_1$ has multiplicity $2$ in $\chi''$, yielding $p+1$ central cyclic coverings of $\M$ for each prime $p$ not dividing $|G|$: three are branched over just two from the three sets of vertices, edges and faces, and $p-2$ are branched over all three sets.

\end{document}